# Algorithme d'aides à la décision pour Optimiser l'Ordonnancement des tâches de maintenance en temps-réel.

Mr Aboussalim Aissam, Mr Mediouni Mohamed
LGII- Laboratoire du Génie Industriel et informatique Ensa Agadir
Email : aboussalimaissam@hotmail.fr

**Résumé :**
L'approche présentée dans cet article est une contribution à la recherche portant sur l'optimisation de l'ordonnancement des tâches.
Ce travail répond aux objectifs suivants : Décision en temps-réel (planification et ordonnancement), Allocation/affectation de ressources, Assurer une disponibilité minimum requise, Minimiser les temps retard en insérant les tâches dynamiques « maintenance corrective ».

**MOT-CLÉS :** tâches, ordonnancement, optimisation, ressource, coût

## 1. INTRODUCTION

Le problème d'ordonnancement consiste à organiser dans le temps la réalisation des tâches, compte tenu de contraintes temporelles et de contraintes portant sur l'utilisation et la disponibilité de ressources requises. L'optimisation de l'ordonnancement des tâches est devenue un domaine essentiel de la recherche. Il contribue à aider les entreprises à rentabiliser aux mieux leurs systèmes de production.

Nous nous intéressons plus particulièrement aux stratégies d'insertion des taches dynamiques pour les quelles les décisions (entretenir, réparer, remplacer, inspecter, …) sont prises sur la base d'informations de surveillance en ligne.

La modélisation des systèmes à composants multiples avec dépendances et sur les techniques d'optimisation multicritères. Dans la plupart des articles dédiés à ce type de problème, les auteurs supposent qu'il n'y a pas de temps mort entre deux tâches consécutives. De plus, les problèmes traités se limitent à l'ordonnancement de tâches sur une seule machine. Tout d'abord de présenter un algorithme qui garantisse le respect des contraintes temporelles et des ressources et assure l'ordonnancement des tâches ainsi que l'affectation des ressources.

Afin de réduire les coûts de maintenance, les travaux de Calibra *et al*. [1] ont traités des problèmes d'allocation de fiabilité, L'amélioration des coûts de maintenance et de la disponibilité peut aussi être obtenue par la gestion optimale des ressources humaines, de même que Parmi les rares travaux qui abordent le problème typique de planification optimale des ressources et activités de maintenance, il y a notamment les travaux de Graves et Lee [2], Lee et Chen [3], Ce dernier groupe de travaux est plutôt fondé sur des approches de gestion de production. Ainsi Les travaux par Cervin d'ordonnancement et planification des taches en temps réelles [4],

L'approche présentée dans ce travail a pour objectif d'affecter les ressources aux tâches de façon à minimiser, sur un horizon *H, le* coût des retards ou déviations, de proposer un algorithme permettant de trier, insérer les tâches dynamiques en temps réel, offrir une interface utilisateur affin d'ordonnancer les tâches en temps réel et permettre une utilisation maximale des ressources pour un coût perdu minimal.

Ce papier est organisé de la manière suivante, dans la deuxième section, on présente la position du problème,, la section 3, on donne la solution proposée, et dans la section 4 on présente les résultats expérimentaux et enfin conclusion.

## 2. Position du problème

Le problème d'ordonnancement de la maintenance ou de la production est NP-difficile.il est classé en 2categories, ordonnancement en ligne et ordonnancement hors-ligne. Un ordonnancement hors-ligne signifie que la séquence d'ordonnancement est prédéterminée à l'avance : dates de début d'exécution des tâches. En pratique, l'ordonnancement prend la forme d'un plan hors-ligne (ou statique), exécuté de façon répétitive. L'ordonnanceur s'appelle dans ce cas un *séquenceur*. Un ordonnancement en-ligne correspond au déroulement d'un algorithme qui tient compte des tâches effectivement présentes dans la file d'ordonnancement. lors de chaque décision d'ordonnancement

On doit minimiser les dates d'achèvement de la réalisation des tâches, optimiser l'utilisation des ressources sous des contraintes de durées. On plus il faut insérer des taches dynamiques en temps réel. Pour modéliser ce problème par des algorithmes considérés comme NP difficile on a procédé à donner des heuristiques qui permettent d'obtenir des solutions approchées par rapport à la solution optimale

Au niveau de l'organisation du système, nous adoptons la structure suivante : un centre de compétence, qui évalue les ressources humaines et les classe par bilan de compétence pour une période déterminée.

Ce modèle est une suite de l'article d'Algorithme d'arrête gauche [5] adapté pour l'ordonnancement des tâches de maintenance. Des fenêtres optimales permettent d'affecter la tâche courante dans la meilleure fenêtre élémentaire et de calculer à chaque fenêtre élémentaire le coût perdu. La valeur ajoutée de notre approche est d'affecter les ressources les plus compétentes aux premiers tâches et ceci pour achever les tâches dans les meilleurs délais et minimiser la somme pondérée des temps mort, en insérant les tâches dynamiques dans les fenêtres avoisines de leur durées, ainsi l'indisponibilité de chaque site et le coût perdu doit être minimal et de programmer un nombre illimité de taches

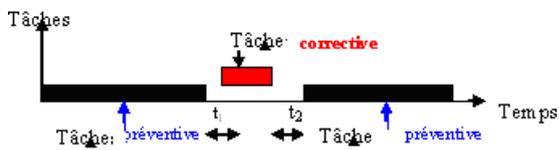

**FIG1 : Insertion de la tache dynamique dans la fenêtre élémentaire**

Note but est d'insérer une tâche dynamique telle que lim t1 ⟶ 0 et lim t2 ⟶ 0

### a- Hypothèses

Pour formuler le problème, nous avons adopté les hypothèses suivantes :

- Nous nous positionnons dans le cas des tâches préventives pour planifier les tâches sur un horizon H.
- Les équipes ont des niveaux de compétences différents classés par un service d'examens responsables des bilans de compétence.
- On procède à appliquer notre modèle sur une seule machine
- Les tâches ont des priorités égales.
- Les dates d'intervention et les durées de tâches sont calculées sur la base de l'analyse statistique des pannes.
- Chaque tâche achevée doit être supprimée.
- Les dates d'échéance de tâches sont réparties dans l'ordonnancement optimal, tous les temps libres entre les tâches seront remplit par des tâches courantes pour minimiser Le coûts global

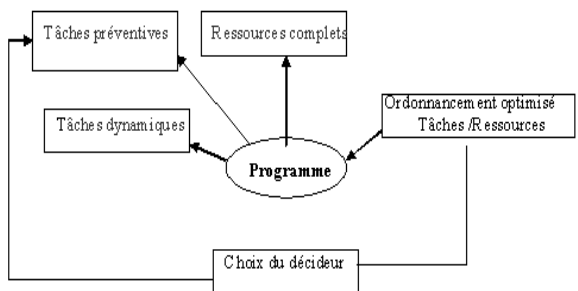

**fiG2 : modèle d'affectation et d'ordonnancement dynamiques d'une nouvelle tâche**

### b- Formulation des hypothèses

L'ordonnancement d'un ensemble de n tâches periodiques est faisable

$$\tau = \{\tau_1, \ldots \tau_i, \ldots \tau_n\}$$

**Si le facteur d'utilisation**

$$U = \sum_{i=1}^{n} C_i / P_i \quad [6]$$

Vérifie la relation suivante

$$U \leq n \cdot (\sqrt[n]{2} - 1) \quad [6]$$

une condition suffisante d'ordonnançabilité d'après l'algorithme RM [6]

$C_i$ : Durée d'exécution maximale ; $P_i$ : période de la tâche

Le service est composé de m ressources humaines (i=1…m), chacune Caractérisée par un profil de compétence. Chaque ressource a donc un niveau de qualification correspondant à chaque type de tâche. C'est ce niveau de qualification qui fait que, d'une ressource à une autre, le temps de traitement ne sera pas le même.

Où $C_{ik}$ est le taux de compétence de la ressource i dans la compétence requise pour le traitement de la tâche du type k. Il est donc possible de représenter cela par une matrice dans laquelle, pour chaque type de tâche et pour chaque ressource, on obtient le taux de Compétence correspondant.. Cette partie présente la formulation du problème d'affectation qui peut être développée dans le cas où l'on impose de satisfaire, pour un processus donné, la contrainte Performance_ Processus ≤ OB (durée, coût). [7]

OB :l'objectif de performance imposé au processus Elle peut être adaptée sans difficulté au cas où la contrainte définie est de type

Performance_ Processus ≥ OB (qualité). [7]

$$\begin{bmatrix} C_{1,1} & \cdots & C_{1,o} \\ \vdots & \ddots & \vdots \\ C_{m,1} & \cdots & C_{m,o} \end{bmatrix} \quad [8]$$

Le coût des activités est la somme des coûts de toutes les tâches et les coûts perdus entre taches qu'on peut appeler fenêtre. Le coût perdu est le produit du temps de déplacement et le coût d'une heure de déplacement :

$C_{k\ perdu} = [T_{K+1} - T_K] \times C_{heure}$

$T_K$ : Temps d'exécution de la tâche k

La fonction globale s'écrit comme :

$$F_{OG} = \sum_i C_{Ti} + \sum_k C_{k\ perdu}$$

CTi : cout d'une tache i

$$F_{opti\_tâches} = \sum_i C_{Ti} = \sum_{l=1} W_i \theta_i + h_i \theta_{i-} + C_{i0} \quad (9)$$

« Cette fonction décide pour chaque tâche de dévier positivement ou négativement respectivement par rapport à sa date due au plus tôt ou sa date due au plus tard, soit de co

Chaque tâche i se caractérise par une durée opératoire notée Ti

$\theta_{i-}$ :: désigne la déviation négative de la tâche i par rapport à sa date due au plus tôt ou également l'avance.

$\theta_{i+}$ :: désigne la déviation positive de la tâche i par rapport à sa date due au plus tard ou également le retard.

Une pénalité de retard Wi, une pénalité d'avance hi, Ci0 le coût minimal de maintenance préventive

$$F_{OG} = \sum_{l=1} W_i \theta_i + h_i \theta_{i-} + C_{i0} + \min \sum_K [T_{K+1} - T_K] \times C_{heure}$$

Le coût des activités est la somme des coûts des tâches ajouté aux couts perdus entre taches k+1et k

### 3. Solution proposée

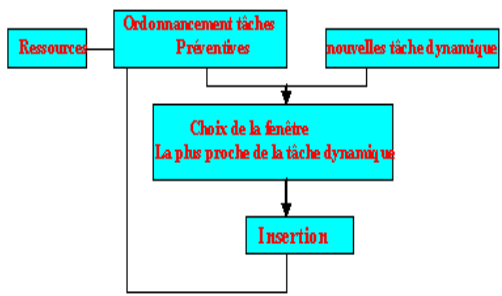

**Fig3 Schéma de principe optimisation d'ordonnancement des tâches**

La solution proposée reçoit en entrée une liste L contenant des tâches avec leurs caractéristiques (ri, di, ci).

  ri : date début de tâche (ri).
  di : date fin de tâche
  ci est le coût de la tâche :

- ✓ Optimisation d'une seule tâche afin d'être planifiée dans son intervalle optimal
- ✓ Tri des tâches par ordre croissant par date début et par temps d'exécution (di-ri)
- ✓ Affectation des ressources les plus compétentes en temps d'exécution (di-ri) le plus grand.
- ✓ Calcul des coûts des tâches et de coûts perdu par la fonction FOG
- ✓ Calculer temps d'exécution de la tâche dynamique (dj-rj)
- ✓ Chercher et calculer le temps de toutes les fenêtres« temps perdu entre les taches »
- ✓ Insertion des tâches dynamique dans la meilleure fenêtre dont son temps d'exécution est approximativement égale.
- ✓ affecter les ressources les plus compétentes disponibles aux premières tâches dynamiques les plus longues.

**EXEMPLE D'APPLICATION**
- Insertion des nouvelles taches courantes non prévues ou en retards
- Robot d'exploration : Ce rebot doit se déplacer dans un environnement non connu (zone radioactive, planète, épave sous la mer).

- Système de contrôle arien : ce système rassemble l'information sur l'état de chaque avion via un ou plusieurs radars actifs. Ce radar interroge chaque avion périodiquement. Le système traite les messages reçus de l'avion et stocke les informations dans une base de données. Ces informations sont prises en compte et traité par des processeurs graphiques. Au même temps, un système de surveillance analyse de façon continue le scenario et alerte les opérateurs toutefois qu'il détecte une collision possible.

Tous ces exemples citées on doit avoir des temps libre entre taches automatisées afin d'intervenir

## Algorithme

**Affectation des tâches (Ti)**  // Ti : Tâche i
**Début**
   Si (P=∅)    alors    // P : Pile des Taches
          Pile ← Ti   // début de la pile
   Sinon    Aller_Fin_Pile (Pile) P ← Ti
          // Ti est ajouté à la fin de la pile
    **Fin**Si
**FIN**

**Affecter_Pile_Ressource (Ri)**  // Ri : ressource i
**Début**
   Si (R⟨⟩ ∅)   alors  Pile ← Ri
        // R : Pile des   ressources
         sinon
         Aller_Fin_Ressource (R)
           R ← Ri
    **FINSI**
**Fin**

**Tri_Ressource (R)**   //R : Pile des ressources, une ressource se caractérise par une note et une variable disponible
**Début**
       Ri ← Premiere_ressource
        Tant que (Ri.disponible = False)
           Ri ← Ressource_suivante

         **Fin**Tque
     **Si** (Ri .disponible = vrai ) alors
       Affecter_Ressource_Tache  (Ri , P)
           **Sinon**
           Ecrire «Ressource non  Disponible »
      **Fin Si**
  F**in**

**Affecter_Ressource_Tache (Ri, pile P)**
//Ri : Ressource la plus Compétente
disponibleVariable  $P_gT$ : Tâche   // $P_gT$  la plus grande tâche
  Variable h ,pos,N   //h compteur de la pile P
   // N nombre de tâches
// pos est la position de tâche dans la pile di et ri sont les paramètres de la tâche  (début et fin)

```
Trier_tâche ( )
  h ← 2
  pos ← 1
  P_gT ← T_1.di - T_1.ri
DEBUT
  Tant que (h ≤ N)
    Si ((T_h.d_h - T_h.r_h) > P_gT) alors
      P_gT ← T_h
      Pos ← h
    Finsi
    h ← h+1
  FinTque
  Ri.disponible ← False
Fin
```
// Disponible est un paramètres qui vaut true si la ressource n'est pas prise sinon false

**Calcul _cout_perdu ( )** //Temps d'arrêt
**Début**
  Somme=0
  Pour i=1 à N faire
    Somme=somme + $T_{i+1}.d - T_i.r$
  Fin Pour
  // Afficher « coût perdu = », somme*Ci0
  **Calcul _coût_ perdu** ← somme*Ci0
**Fin**
  F_avant ← **Calcul _coût_ perdu( )**

**Apres Optimisation  //insertion des tâches dynamiques**
Inser_Tache_D  (Td ; P)
Début
  Var X, T : tâche
  Var : Durée _TD, pos : ENTIER
  Durée –Td ← Td.di - Td.ri
  T ← première_ tache (Pile)
  Pos ← 1
  i ← 1
Tant que (T < > Fin_ pile)
  X ← Tache _ Suivante (Pile)
  Min ← X.ri – T.di
  SI (Min >= Durée _ Td) alors
    Pos ← i
    T ← X
    X ← Tache _suivante (pile)
    SI (Min >= (X.ri – T.di)) et
       (X.ri – T.di) > Durée- TD)
      Min ← X.ri – T.di
      T ← X
      Pos ← i
    fsi
  fsi
fin
// Se positionner avant la tâche ou il faut insérer
  T ← premier Tache (Rle)
  Pour i = 1 à pos
    T ← Tache_ suivant (Rle)
  Fpour
  Z ← T.suivant
  T.Suivant ← Td
  Td.suivant ← Z.suivant

**Calcul _cout_perdu ( )**
  //Temps d'arrêt  //après optimisation
**Début**
  Somme=0
  Pour i=1 a N faire
    Somme=somme + $T_{i+1}.d - T_i.r$
  Fin Pour
  // Afficher « coût perdu = », somme*Ci0
  **Calcul _coût_ perdu** ← somme*Ci0
**Fin**
  F_apres ← **Calcul _coût_ perdu( )**
  //après insertion des tâches dynamiques
  // **Tg** ← F_avant - F_apres

## 4. RÉSULTATS

L'algorithme d'affectation des tâches a été programmé par JAVA.et la disposition des tâches ressources par ms-Project. Pour évaluer ce programme, nous avons considéré un ensemble de 10 tâches préventives à ordonnancer en fonction du temps d'exécution et des ressources. Puis on a inséré dans un premier temps 3 tâches dynamiques et dans un deuxième temps 9 tâches dynamiques.

### a- Taches avant modélisation

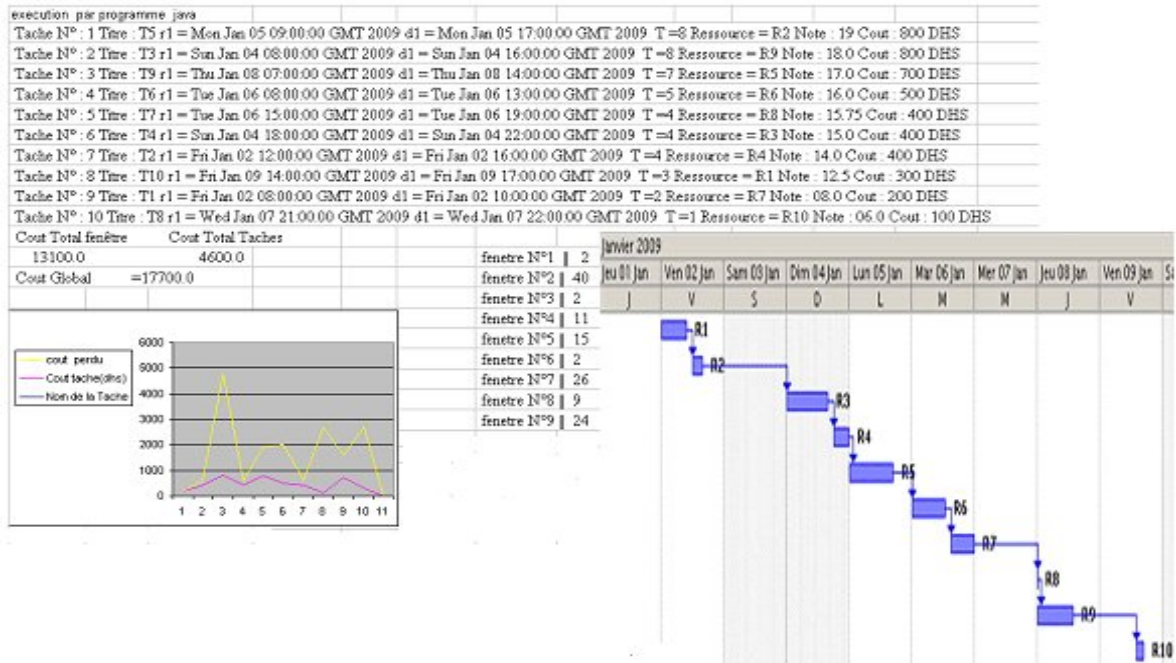

**Constatation1**

Les tâches sont triées par ordre croissant du temps d'exécution .la tâche la plus longue est. Affectée à la ressource la plus compétente .le programme calcule le coût de chaque tâche et le coût perdu entre tâches

### b- Taches après modélisation(optimisation) INSERTION DE 3TACHES DYNAMIQUES

Tache N° :1 Titre : T1 r1 = Fri Jan 02 08:00:00 GMT 2009 d1 = Fri Jan 02 10:00:00 GMT 2009  T =2 Resource = R7 Note : 08.0 Coot : 200 DHS
Tache N° :2 Titre : TD1 r1 = Fri Jan 02 10:00:00 GMT 2009 d1 = Fri Jan 02 11:00:00 GMT 2009  T =1 Ressource = R2 Note : 19 Coot: 100 DHS
Tache N° :3 Titre : T2 r1 = Fri Jan 02 12:00:00 GMT 2009 d1 = Fri Jan 02 16:00:00 GMT 2009  T =4 Resource = R4 Note : 14.0 Coot : 400 DHS
Tache N° :4 Titre : TD2 r1 = Fri Jan 02 17:00:00 GMT 2009 d1 = Sun Jan 04 08:00:00 GMT 2009  T =39 Resource = R9 Note : 18.0 Cout : 3900
Tache N° :5 Titre : T3 r1 = Sun Jan 04 08:00:00 GMT 2009 d1 = Sun Jan 04 16:00:00 GMT 2009  T =8 Resource = R2 Note : 19 Cout : 800 DHS
Tache N° :6 Titre : T4 r1 = Sun Jan 04 18:00:00 GMT 2009 d1 = Sun Jan 04 22:00:00 GMT 2009  T =4 Resource = R3 Note : 15.0 Cout : 400 DHS
Tache N° :7 Titre :TD3 r1 = Sun Jan 04 23:00:00 GMT 2009 d1 = Mon Jan 05 07:00:00 GMT 2009  T =8 Resource = R5 Note : 17.0 Coot  800 DHS
Tache N° :8 Titre : T5 r1 = Mon Jan 05 08:00:00 GMT 2009 d1 = Mon Jan 05 16:00:00 GMT 2009  T =8 Resource = R9 Note : 18.0 Cout : 800 DHS
Tache N° :9 Titre : T6 r1 = Tue Jan 06 08:00:00 GMT 2009 d1 = Tue Jan 06 13:00:00 GMT 2009  T =5 Resource = R6 Note : 16.0 Cout : 500 DHS
Tache N° :10 Titre : T7 r1 = Tue Jan 06 15:00:00 GMT 2009 d1 = Tue Jan 06 19:00:00 GMT 2009  T =4 Resource = R8 Note : 15.75 Coot : 400 DHS
Tache N° :11 Titre : T8 r1 = Wed Jan 07 21:00:00 GMT 2009 d1 = Wed Jan 07 22:00:00 GMT 2009  T =1 Resource = R10 Note : 06.0 Coot : 100 DHS
Tache N° :12 Titre : T9 r1 = Thu Jan 08 07:00:00 GMT 2009 d1 = Thu Jan 08 14:00:00 GMT 2009  T =7 Ressource = R5 Note :17.0 Coot : 700 DHS
Tache N° :13 Titre : T10 r1 = Fri Jan 09 14:00:00 GMT 2009 d1 = Fri Jan 09 17:00:00 GMT 2009  T =3 Resource = R1 Note : 12.5 Coot : 300 DHS

Cout Total fen?tre    Cout Total Taches
8300.0                9000.0
Co?t Global    =17300.0

fenetre N°1 || 0
fenetre N°2 || 1
fenetre N°3 || 1
fenetre N°4 || 0
fenetre N°5 || 2
fenetre N°6 || 1
fenetre N°7 || 1
fenetre N°8 || 16
fenetre N°9 || 2
fenetre N°10 || 26
fenetre N°11 || 9
fenetre N°12 || 24

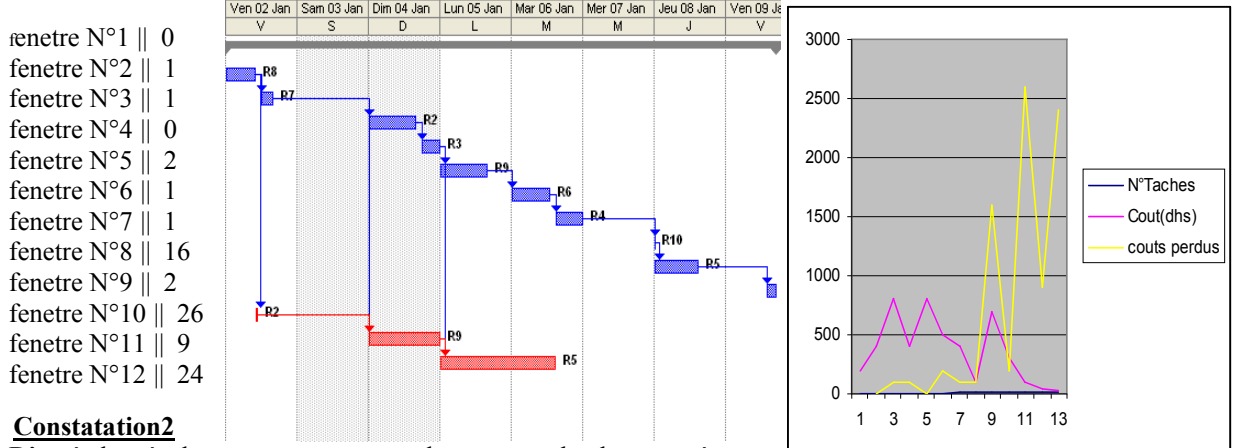

**Constatation2**

D'après les résultats, on remarque que la ressource la plus compétente est affectée a la tache préventive la plus longue et c'est cette ressource qui se réaffecte à la tâche dynamique la plus longue. On remarque d'après la courbe que le coût perdu après optimisation a été atténué de **37%.**

### c- Taches après modélisation(optimisation) INSERTION DE 9TACHES DYNAMIQUES

Tache N° :1 Titre : T1 r1 =Fri Jan 02 08:00:00 GMT 2009 d1= Fri Jan 02 10:00:00 GMT 2009T =2 Resource = R7 Note : 08.0 Coot :200DHS
Tache N° :2 Titre : TD1 r1 = Fri Jan 02 10:00:00 GMT 2009 d1 = Fri Jan 02 11:00:00 GMT 2009  T =1 Resource = R1 Note : 12.5 Coot : 300 DHS
Tache N° :3 Titre : T2 r1 = Fri Jan 02 12:00:00 GMT 2009 d1 = Fri Jan 02 16:00:00 GMT 2009  T =4 Resource = R4 Note : 14.0 Coot : 400 DHS
Tache N° :4 Titre : TD2 r1 = Fri Jan 02 17:00:00 GMT 2009 d1 = Sun Jan 04 08:00:00 GMT 2009  T =39 Resource = R2 Note : 19 Cout : 3900 DHS
Tache N° :5 Titre : T3 r1 = Sun Jan 04 08:00:00 GMT 2009 d1 = Sun Jan 04 16:00:00 GMT 2009  T =8 Ressource = R2 Note : 19 Cout : 800 DHS
Tache N° :6 Titre : TD4 r1 = Sun Jan 04 16:00:00 GMT 2009 d1 = Sun Jan 04 17:00:00 GMT 2009  T =1 Resource = R7 Note : 08.0 Coot : 200 DHS
Tache N° :7 Titre : T4 r1 = Sun Jan 04 18:00:00 GMT 2009 d1 = Sun Jan 04 22:00:00 GMT 2009  T =4 Resource = R3 Note : 15.0 Coot : 400 DHS
Tache N° :8 Titre : TD3 r1 = Sun Jan 04 23:00:00 GMT 2009 d1 = Mon Jan 05 07:00:00 GMT 2009  T =8 7Ressource = R5 Note : 17

Tache N° :9 Titre : T5 r1 = Mon Jan 05 08:00:00 GMT 2009 d1 = Mon Jan 05 16:00:00 GMT 2009  T =8 Ressource = R9 Note : 18.0 Cout : 800 DHS
Tache N° :10 Titre : TD5 r1 = Mon Jan 05 16:00:00 GMT 2009 d1 = Mon Jan 05 18:00:00 GMT 2009  T =2 Ressource = R3 Note : 15.0 Coot : 200 DHS
Tache N° :11 Titre : T6 r1 = Tue Jan 06 08:00:00 GMT 2009 d1 = Tue Jan 06 13:00:00 GMT 2009  T =5 Ressource = R6 Note : 16.0 Cout : 500 DHS
Tache N° :12 Titre : T7 r1 = Tue Jan 06 15:00:00 GMT 2009 d1 = Tue Jan 06 19:00:00 GMT 2009  T =4 Ressource = R8 Note : 15.75 Coot : 400 DHS
Tache N° :13 Titre : TD6 r1 = Tue Jan 06 20:00:00 GMT 2009 d1 = Tue Jan 06 23:00:00 GMT 2009  T =3 Ressource = R8 Note : 15.75 Coot : 300 DHS
Tache N° :14 Titre : TD7 r1 = Wed Jan 07 08:00:00 GMT 2009 d1 = Wed Jan 07 17:00:00 GMT 2009  T =9 Ressource = R9 Note : 18.0 Cout : 900 DHS
Tache N° :15 Titre : T8 r1 = Wed Jan 07 21:00:00 GMT 2009 d1 = Wed Jan 07 22:00:00 GMT 2009  T =1 Ressource = R10 Note : 06.0 Coot : 100 DHS
Tache N° :16 Titre : T9 r1 = Thu Jan 08 07:00:00 GMT 2009 d1 = Thu Jan 08 14:00:00 GMT 2009  T =7Ressource = R5 Note : 17.0 Coot : 700 DHS
Tache N° :17 Titre : TD8 r1 = Thu Jan 08 15:00:00 GMT 2009 d1 = Thu Jan 08 23:00:00 GMT 2009  T =8 Ressource = R6 Note : 16.0 Cout : 800 DHS
Tache N° :18 Titre : T10 r1 = Fri Jan 09 14:00:00 GMT 2009 d1 = Fri Jan 09 17:00:00 GMT 2009  T =3 Ressource = R1 Note : 12.5 Coot : 300 DHS
Tache N° :19 Titre : TD9 r1 = Fri Jan 09 18:00:00 GMT 2009 d1 = Fri Jan 09 20:00:00 GMT 2009  T =2 Ressource = R4 Note : 14.0 Coot : 100 DHS

Cout Total fen?tre    Cout Total Taches
6100.0                10000.0

## Co?t Global     =16100.0

fenetre N°1 || 0      fenetre N°7 || 1     fenetre N°13 || 9
fenetre N°2 || 1      fenetre N°8 || 1     fenetre N°14 || 4
fenetre N°3 || 1      fenetre N°9 || 0     fenetre N°15 || 9
fenetre N°4 || 0      fenetre N°10 || 14   fenetre N°16 || 1
fenetre N°5 || 0      fenetre N°11 || 2    fenetre N°17 || 15
fenetre N°6 || 1      fenetre N°12 || 1    fenetre N°18 || 1

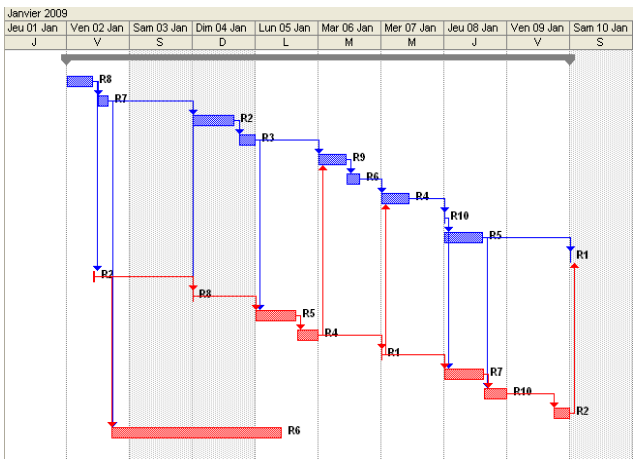 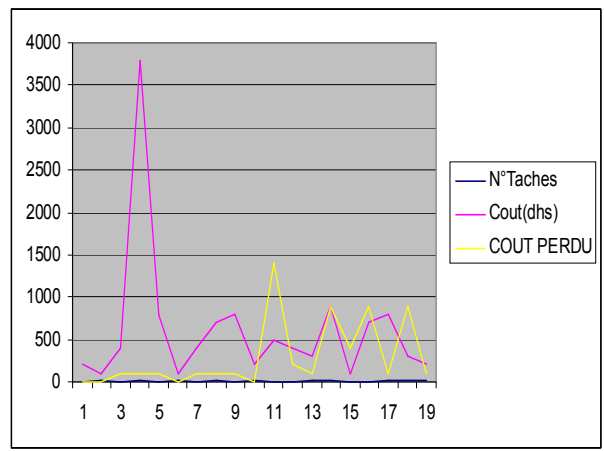

**Constatation3**

D'après les résultats, on remarque que la ressource la plus compétente est affecté à la tache préventive la plus longue et c'est cette ressource qui se réaffecte a la tache dynamique la plus longue. On remarque d'après la courbe que le cout perdu après optimisation a été atténué de **54%**

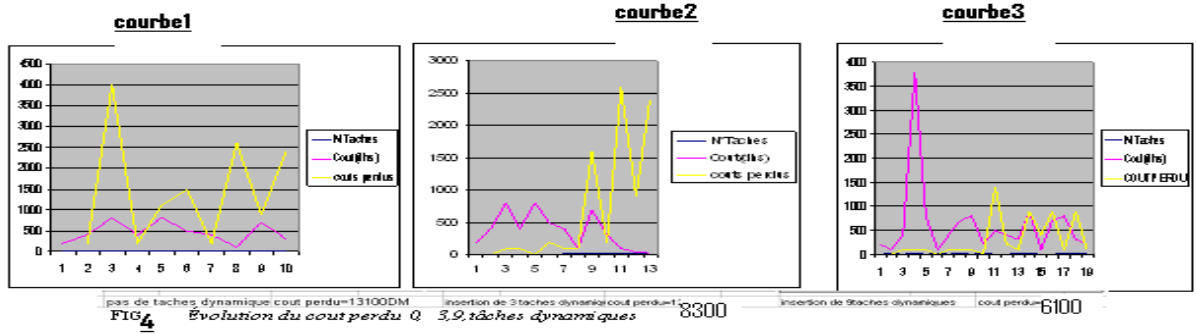

FIG 4 Évolution du cout perdu Q 3,9 tâches dynamiques

**Interprétation des courbes**

En comparant les 3courbes on constate que le cout perdu chute en insérant de plus en plus de taches dynamiques
Cout perdu après optimisation (9taches dynamiques) =6100DHS< Cout perdu après optimisation(3 taches dynamiques) =8300DHS < Cout perdu avant optimisation=13100 DHS
Courbe 1 : coût des taches préventives et cout perdu entre taches
Courbe 2 : Le coût perdu est quasiment nul jusqu' à la tâche préventive  N°7 commence à augmenter ce qui explique il faut programmer des tâches dynamiques à partir de cette  tâche
Courbe 3 : Le coût perdu est quasiment nul jusqu' à la tache préventive  N°10 ce qui explique que le programme a inséré des taches dynamiques .au dessous de cette tâche

## 4. **CONCLUSION**

Ce travail traite un problème d'optimisation d'ordonnancement des tâches et d'affectation des ressources, qui a pour objectif la réduction du coût perdu. Nous avons présenté un algorithme qui peut ordonnancer n tâches préventives, (n-1) tâches dynamiques. Cette méthode a montré de bons résultats d'optimisation du coût perdu en insérant de plus en plus de tâches courantes.

# Annexes

**Tableau1** — saisie des durees des taches preventives et calcul des fenetres

| N°Taches | Durée (H) | Début | Fin | Cout(dhs) | Nom Ressourc | type taches | couts perdus | fenetres |
|---|---|---|---|---|---|---|---|---|
| 1 | 2 | 2/1/09 8:00 | 2/1/09 10:00 | 200 | R8=15,75 | preventive | | |
| 2 | 4 | 2/1/09 12:00 | 2/1/09 16:00 | 400 | R7=8 | preventive | 200 | fenetre N°1 |
| 3 | 8 | 4/1/09 8:00 | 4/1/09 16:00 | 800 | R2=19 | preventive | 4000 | fenetre N°2 |
| 4 | 4 | 4/1/09 18:00 | 4/1/09 22:00 | 400 | R3=15 | preventive | 200 | fenetre N°3 |
| 5 | 8 | 5/1/09 8:00 | 5/1/09 16:00 | 800 | R9=18 | preventive | 1100 | fenetre N°4 |
| 6 | 5 | 6/1/09 8:00 | 6/1/09 13:00 | 500 | R6=16 | preventive | 1500 | fenetre N°5 |
| 7 | 4 | 6/1/09 15:00 | 6/1/09 19:00 | 400 | R4=14 | preventive | 200 | fenetre N°6 |
| 8 | 1 | 7/1/09 21:00 | 7/1/09 22:00 | 100 | R10=6 | preventive | 2600 | fenetre N°7 |
| 9 | 7 | 8/1/09 7:00 | 8/1/09 14:00 | 700 | R5=17 | preventive | 900 | fenetre N°8 |
| 10 | 3 | 9/1/09 14:00 | 9/1/09 17:00 | 300 | R1=12,5 | preventive | 2400 | fenetre N°9 |
| cout perdu(avant optimisation)=131*100 =13100DHS | | | | | | | 13100 | |

**Tableau2** — insertion de 3 taches dynamiques

| N°Taches | Durée (H) | Début | Fin | Cout(dhs) | Nom Ressourc | type taches | couts perdus | fenetres |
|---|---|---|---|---|---|---|---|---|
| 1 | 2 | 2/1/09 8:00 | 2/1/09 10:00 | 200 | R8=15,75 | preventive | | |
| 11 | 1 | 2/1/09 10:00 | 2/1/09 11:00 | 100 | R2=19 | dynamique | 0 | fenetre N°1 |
| 2 | 4 | 2/1/09 12:00 | 2/1/09 16:00 | 400 | R7=8 | preventive | 100 | fenetre N°2 |
| 3 | 8 | 4/1/09 8:00 | 4/1/09 16:00 | 800 | R2=19 | preventive | 100 | fenetre N°3 |
| 12 | 15 | 4/1/09 8:00 | 4/1/09 23:00 | 1500 | R9=18 | dynamique | 0 | fenetre N°4 |
| 4 | 4 | 4/1/09 18:00 | 4/1/09 22:00 | 400 | R3=15 | preventive | 200 | fenetre N°5 |
| 13 | 22 | 4/1/09 23:30 | 5/1/09 21:30 | 2200 | R5=17 | dynamique | 100 | fenetre N°6 |
| 5 | 8 | 5/1/09 8:00 | 5/1/09 16:00 | 800 | R9=18 | preventive | 100 | fenetre N°7 |
| 6 | 5 | 6/1/09 8:00 | 6/1/09 13:00 | 500 | R6=16 | preventive | 1600 | fenetre N°8 |
| 7 | 4 | 6/1/09 15:00 | 6/1/09 19:00 | 400 | R4=14 | preventive | 200 | fenetre N°9 |
| 8 | 1 | 7/1/09 21:00 | 7/1/09 22:00 | 100 | R10=6 | preventive | 2600 | fenetre N°10 |
| 9 | 7 | 8/1/09 7:00 | 8/1/09 14:00 | 700 | R5=17 | preventive | 900 | fenetre N°11 |
| 10 | 3 | 9/1/09 14:00 | 9/1/09 17:00 | 300 | R1=12,5 | preventive | 2400 | fenetre N°12 |
| | | | | | | | 8300 | |

**Tableau3** — insertion de 9 taches dynamiques

| N°Taches | Durée (H) | Début | Fin | Cout(dhs) | Nom Ressourc | fenetres | COUT PERDU | fenetres |
|---|---|---|---|---|---|---|---|---|
| 1 | 2 | 2/1/09 8:00 | 2/1/09 10:00 | 200 | R8=15,75 | 0 | 0 | |
| 11 | 1 | 02/01/2009 10:00 | 02/01/2009 11:00 | 100 | R9=18 | 0 | 0 | fenetre N°1 |
| 2 | 4 | 2/1/09 12:00 | 2/1/09 16:00 | 400 | R7=8 | 1 | 100 | fenetre N°2 |
| 19 | 38 | 02/01/2009 17:00 | 04/01/2009 07:00 | 3800 | R6=16 | 1 | 100 | fenetre N°3 |
| 3 | 8 | 4/1/09 8:00 | 4/1/09 16:00 | 800 | R2=19 | 0 | 0 | fenetre N°4 |
| 12 | 1 | 04/01/2009 16:00 | 04/01/2009 17:00 | 100 | R8=15,75 | 0 | 0 | fenetre N°5 |
| 4 | 4 | 4/1/09 18:00 | 4/1/09 22:00 | 400 | R3=15 | 1 | 100 | fenetre N°6 |
| 13 | 7 | 04/01/2009 23:00 | 05/01/2009 07:00 | 700 | R3=15 | 1 | 100 | fenetre N°7 |
| 5 | 8 | 5/1/09 8:00 | 5/1/09 16:00 | 800 | R9=18 | 1 | 100 | fenetre N°8 |
| 14 | 2 | 05/01/2009 16:00 | 05/01/2009 18:00 | 200 | R4=14 | 0 | 0 | fenetre N°9 |
| 6 | 5 | 6/1/09 8:00 | 6/1/09 13:00 | 500 | R6=16 | 14 | 1400 | fenetre N°10 |
| 7 | 4 | 6/1/09 15:00 | 6/1/09 19:00 | 400 | R4=14 | 2 | 200 | fenetre N°11 |
| 15 | 3 | 06/01/2009 20:00 | 06/01/2009 23:00 | 300 | R1=12,5 | 1 | 100 | fenetre N°12 |
| 16 | 9 | 07/01/2009 08:00 | 07/01/2009 17:00 | 900 | R7=8 | 9 | 900 | fenetre N°13 |
| 8 | 1 | 7/1/09 21:00 | 7/1/09 22:00 | 100 | R10=6 | 4 | 400 | fenetre N°14 |
| 9 | 7 | 8/1/09 7:00 | 8/1/09 14:00 | 700 | R9=18 | 9 | 900 | fenetre N°15 |
| 17 | 8 | 08/01/2009 15:00 | 08/01/2009 23:00 | 800 | R10=6 | 1 | 100 | fenetre N°16 |
| 10 | 3 | 9/1/09 14:00 | 9/1/09 17:00 | 300 | R1=12,5 | 15 | 1500 | fenetre N°17 |
| 18 | 2 | 09/01/2009 18:00 | 09/01/2009 20:00 | 200 | R2=19 | 1 | 100 | fenetre N°18 |
| | | | | | | | 6100 | |